\documentclass[11pt]{article}

\renewcommand{\a }{\alpha }

\renewcommand{\l}{\lambda }
\newcommand{\media}{\m\kern12mu\hbox{\vrule height9pt depth-3.2pt width5pt} \m\kern-20mu\int}

\newcommand{\HH }{\mathrm{H}}

\newcommand{\half}{\frac{1}{2}}

\renewcommand{\l }{\lambda }

\newcommand{\m }{\mu }

\newcommand{\bbR }{\mathbb{R}}

\usepackage{amsmath,amsthm,amsfonts,amssymb,latexsym}
\usepackage{amsmath}

\usepackage{a4wide}

\newtheorem{thm}{Theorem}[section]

\newtheorem{prop}[thm]{Proposition}
\newtheorem{ex}[thm]{Example}

\newtheorem{theorem}{Theorem}[section]

\newtheorem{lemma}[theorem]{Lemma}

\newtheorem{proposition}[theorem]{Proposition}

\begin{document}


\date{}\title{The Binary $\mathfrak{aff}(n|1)$-Invariant Differential Operators On Weighted Densities On The Superspace $\mathbb{R}^{1|n}$ And $\mathfrak{aff}(n|1)$-Relative Cohomology}

\author{Khaled Basdouri \thanks{D\'epartement de Math\'ematiques, Facult\'e des Sciences de Gafsa, zarroug 2112,
Tunisie.~~E.mail: khaled.basdouri2016@gmail.com}\and
Salem Omri \thanks{D\'epartement de Math\'ematiques, Facult\'e des Sciences de
Gafsa, zarroug 2112, Tunisie.~~~~E.mail: omri\_salem@yahoo.fr
}\and
Wissal Swilah \thanks{D\'epartement de Math\'ematiques, Facult\'e des Sciences de
Gafsa, zarroug 2112, Tunisie.~~~~E.mail: swilah\_wissal@yahoo.fr
}}
\maketitle
\begin{abstract}
We consider the $\mathfrak{aff}(n|1)-$module structure on the spaces
of differential bilinear operators acting on the superspaces of weighted
densities. We classify $\mathfrak{aff}(n|1)-$invariant binary differential operators acting on the spaces of weighted densities. This result allows us to compute the first $\mathfrak{aff}(n|1)-$relative differential
cohomology of $\mathcal{K}(n)$ with coefficients in the superspace of linear differential operators acting on the superspaces of weighted densities.

\end{abstract}

\maketitle {\bf Mathematics Subject Classification} (2010). 17B56, 53D55, 58H15.


\section{Introduction}

Let $\mathfrak{vect}(1)$ be the Lie algebra of polynomial vector
fields on $\mathbb{R}$. Consider the $1-$parameter deformation of
the $\mathfrak{vect}(1)-$action  on $\mathbb{R}[x]$:
\begin{equation*}
L_{X\frac{d}{dx}}^\lambda(f)= Xf'+\lambda X'f,
\end{equation*}
where $X, f\in\mathbb{R}[x]$ and $X':=\frac{dX}{dx}$. Denote by
$\mathcal{F}_\l$ the $\mathfrak{vect}(1)-$module structure on $\bbR[x]$
defined by $L^\lambda$ for a fixed $\lambda$. Geometrically,
$\mathcal{F}_\lambda=\left\{fdx^{\lambda}\mid f\in
\mathbb{R}[x]\right\}$ is the space of polynomial weighted
densities of weight $\l\in\mathbb{R}$. The space $\mathcal{F}_\lambda$ coincides
with the space of vector fields, functions
and differential $1-$forms for $\lambda = -1,\, 0$ and $1$,
respectively.

\medskip

Denote by $\mathrm{D}_{\lambda,\mu}:=\hbox{Hom}_\mathrm{diff}({\mathcal{F}}_\lambda,
{\mathcal{F}}_\mu)$ the $\mathfrak{vect}(1)-$module of linear
differential operators with the natural
$\mathfrak{vect}(1)-$action. Feigin and Fuchs \cite{ff} computed
$\mathrm{H}^1_\mathrm{diff}\left(\mathfrak{vect}(1);
\mathrm{D}_{\lambda,\mu}\right)$, where $\mathrm{H}^*_\mathrm{diff}$
denotes the differential cohomology; that is, only cochains given
by differential operators are considered. They showed that
non-zero cohomology $\mathrm{H}^1_\mathrm{diff}\left(\mathfrak{vect}(1);
\mathrm{D}_{\lambda,\mu}\right)$ only
appear for particular values of weights that we call {\it
resonant} and which satisfy $\mu-\lambda\in\mathbb{N}$.

If we restrict ourselves to the Lie subalgebra of
$\mathfrak{vect}(1)$ generated by $\{\frac{d}{dx}, x\frac{d}{dx}\}$,
isomorphic to $\mathfrak{aff}(1)$, we get a
family of infinite dimensional $\mathfrak{aff}(1)-$modules,
still denoted $\mathcal{F}_\lambda$ and $\mathrm{D}_{\lambda,\mu}.$ In \cite{bn}, Basdouri and Nasri
classified all $\mathfrak{aff}(1)-$invariant binary differential operators on $\mathbb{R}$ acting in the spaces
$\mathcal{F}_\lambda$. This allows us to show, in this paper, that nonzero
cohomology ${\mathrm H}^1_{\rm
diff}\left(\mathfrak{vect}(1),\mathfrak{aff}(1);
\mathrm{D}_{\lambda,\mu}\right)$ only appear for particular values
of weights which satisfy
$\mu-\lambda\in\mathbb{N}$ and give explicit basis of
this cohomology space. This space arises
in the classification of $\mathfrak{aff}(1)-$trivial infinitesimal deformations of the
$\mathfrak{vect}(1)-$module ${\cal
S}_{\mu-\lambda}=\bigoplus_{k=0}^\infty \mathcal{F}_{\mu-\lambda-k}$, the
space of symbols of $\mathrm{D}_{\lambda,\mu}$.

\medskip

In this paper we study also the simplest super analog of this problem. Namely, we consider the superspace
$\mathbb{R}^{1|n}$ equipped with the contact structure determined by
a 1-form $\alpha_n$, and the Lie superalgebra $\mathcal{K}(n)$ of
contact polynomial vector fields on $\mathbb{R}^{1|n}$. We introduce
the $\mathcal{K}(n)-$module $\mathbb{F}_\lambda^n$ of
$\lambda$-densities on $\mathbb{R}^{1|n}$ and the
$\mathcal{K}(n)-$module of linear differential operators,
$\mathbb{D}_{\lambda,\mu}^n
:=\mathrm{Hom}_{\rm{diff}}(\mathbb{F}_{\lambda}^n,\mathbb{F}_{\mu}^n)$,
which are super analogues of the spaces $\mathcal{F}_\lambda$ and
$\mathrm{D}_{\lambda,\mu}$, respectively. The Lie superalgebra
$\mathfrak{aff}(n|1)$, a super analogue of $\mathrm{\frak {aff}}(1)$,
can be realized as a subalgebra of $\mathcal{K}(n)$. We classify all
$\mathfrak{aff}(n|1)-$invariant binary differential operators on
$\mathbb{R}^{1|n}$ acting in the spaces $\mathbb{F}_\lambda^n$ for
$n=1$ and $2$. We use this result to compute $\HH^1_{\rm
diff}\left(\mathcal{K}(n),
\mathfrak{aff}(n|1);\mathbb{D}^n_{\lambda,\mu}\right)$ for $n=1$
and $2$. We show that nonzero cohomology $\HH^1_{\rm
diff}\left(\mathcal{K}(n),
\mathfrak{aff}(n|1);\mathbb{D}^n_{\lambda,\mu}\right)$ only appear
for resonant values of weights that satisfy $\mu-\lambda\in{1\over
2}\mathbb{N}$. Moreover, we give explicit basis of these cohomology
spaces.

\section{Definitions and Notations}

\subsection{The Lie superalgebra of contact vector fields on
$\mathbb{R}^{1|n}$}

Let $\mathbb{R}^{1\mid n}$ be the superspace with coordinates
$(x,~\theta_1,\ldots,\theta_n),$ where
$\theta_1,\,\dots,\,\theta_n$ are the odd variables, equipped with
the standard contact structure given by the following $1$-form:
\begin{equation}
\label{a} \a_n=dx+\sum_{i=1}^n\theta_id\theta_i.
\end{equation}
On the space
$\mathbb{R}[x,\theta]:=\mathbb{R}[x,\theta_1,\dots,\theta_n]$, we
consider the contact bracket
\begin{equation}
\{F,G\}=FG'-F'G-\frac{1}{2}(-1)^{|F|}\sum_{i=1}^n{\eta}_i(F)\cdot
{\eta}_i(G),
\end{equation}where
${\eta}_i=\frac{\partial}{\partial
{\theta_i}}-\theta_i\frac{\partial}{\partial x}$ and $|F|$ is the
parity of $F$. Note that the derivations $\eta_i$ are the
generators of n-extended supersymmetry and generate the kernel of
the form (\ref{a}) as a module over the ring of polynomial
functions. Let $\mathrm{Vect_{Pol}}(\mathbb{R}^{1|n})$ be the
superspace of polynomial vector fields on ${\mathbb{R}}^{1|n}$:
\begin{equation*}\mathrm{Vect_{Pol}}(\mathbb{R}^{1|n})=\left\{F_0\partial_x
+ \sum_{i=1}^n F_i\partial_i \mid ~F_i\in\mathbb{R}[x,\theta]~
\text{ for all } i  \right\},\end{equation*} where
$\partial_i=\frac{\partial}{\partial\theta_i}$ and
$\partial_x=\frac{\partial}{\partial x} $, and consider the
superspace $\mathcal{K}(n)$ of contact polynomial vector fields on
${\mathbb{R}}^{1|n}$. That is, $\mathcal{K}(n)$ is the superspace
of vector fields on $\mathbb{R}^{1|n}$ preserving the distribution
singled out by the $1$-form $\alpha_n$: $$
\mathcal{K}(n)=\big\{X\in\mathrm{Vect_{Pol}}(\mathbb{R}^{1|n})~|~\hbox{there
exists}~F\in {\mathbb{R}}[x,\,\theta]~ \hbox{such
that}~{L}_X(\alpha_n)=F\alpha_n\big\}. $$ The Lie superalgebra
$\mathcal{K}(n)$ is spanned by the fields of the form:
\begin{equation*}
X_F=F\partial_x
-\frac{1}{2}\sum_{i=1}^n(-1)^{|F|}{\eta}_i(F){\eta}_i,\;\text{where
$F\in \mathbb{R}[x,\theta]$.}
\end{equation*}
In particular, we have $\mathcal{K}(0)=\mathfrak{vect}(1)$. Observe that
${L}_{X_F}(\alpha_n)=X_1(F)\alpha_n$. The bracket in
$\mathcal{K}(n)$ can be written as:
\begin{equation*}
[X_F,\,X_G]=X_{\{F,\,G\}}.
\end{equation*}

\subsection{The superalgebra $\mathfrak{aff}(1|1)$}

The Lie algebra $\frak {aff}(1)$ is isomorphic to the Lie
subalgebra of $\mathfrak{vect}(1)$ generated by
\begin{equation*}\left\{\frac{d}{dx},\,x\frac{d}{dx}\right\}.\end{equation*}
Similarly, we consider the affine Lie superalgebra (\cite{BIIO})
\begin{equation*}\mathfrak{aff}(1|1)=\text{Span}(X_1,\,X_{x},\,
X_{\theta}),\end{equation*} where
\begin{equation*}(\mathfrak{aff}(1|1))_0=\text{Span}(X_1,\,X_{x})
\quad\text{and}\quad(\mathfrak{aff}(1|1))_1=\text{Span}(
X_{\theta}).\end{equation*}
The new commutation relations are
$$[X_1 ,X_x] = X_1, \qquad [X_x,X_\theta] = -\frac{1}{2}X_\theta,$$
 $$[X_1,X_\theta] = 0, \qquad [X_\theta,X_\theta] =\frac{1}{2}X_1.$$
 More generally, the affine Lie superalgebra $\frak{aff}(n|1)$ can be realized as a
subalgebra of $\mathcal{K}(n)$:
\begin{equation*}\frak{aff}(n|1)=\text{Span}(X_1,\,X_{x},\,X_{\theta_i},\cdots,\,X_{\theta_i\theta_j},,\cdots)
\quad 1\leq i,j\leq
n.\end{equation*}
The Lie superalgebra $\frak{aff}(n-1|1)$ can be realized  as a
subalgebra of $\frak{aff}(n|1)$:
\begin{equation*}
 \frak{aff}(n-1|1)=\Big\{X_F\in\frak{aff}(n|1)~|~\partial_{n}F=0\Big\}.
\end{equation*}
Note that, for any  $i$ in $\{1,2,\dots,n-1\}$, $\frak{aff}(n-1|1)$ is isomorphic to
\begin{equation*}
 \frak{aff}(n-1|1)_i=\Big\{X_F\in\frak{aff}(n|1)~|~\partial_{i}F=0\Big\}.
\end{equation*}

\subsection{Modules of weighted densities}

We introduce a one-parameter family of modules over the Lie
superalgebra $\mathcal{K}(n)$. As vector spaces all these modules
are isomorphic to ${\mathbb{R}}[x,\,\theta]$, but not as
$\mathcal{K}(n)-$modules.

For every contact polynomial vector field $X_F$, define a
one-parameter family of first-order differential operators on
$\mathbb{R}[x,\theta]$:
\begin{equation}
\label{superaction} \mathbb{L}^{\lambda}_{X_F}=X_F + \lambda F',\quad
\l\in\mathbb{R}.
\end{equation}
We easily check that
\begin{equation}
\label{crochet}
[\mathbb{L}^{\lambda}_{X_F},\mathbb{L}^{\lambda}_{X_G}]=\mathbb{L}^\l_{X_{\{F,G\}}}.
\end{equation}
We thus obtain a one-parameter family of $\mathcal{K}(n)-$modules on
$\mathbb{R}[x,\theta]$ that we denote $\mathbb
{F}^n_\lambda$, the
space of all polynomial weighted densities on $\mathbb{R}^{1|n}$
of weight $\l$ with respect to $\a_n$:
\begin{equation}
\label{densities} \mathbb
{F}^n_\l=\left\{F\a_n^{\l} \mid F
\in\mathbb{R}[x,\theta]\right\}.
\end{equation}
In particular, we have $\mathbb
{F}_{\lambda}^0=\mathcal{F}_\lambda$. Obviously the adjoint
$\mathcal{K}(n)-$module is isomorphic to the space of weighted densities
on $\mathbb{R}^{1|n}$ of weight $-1.$

\subsection{Differential operators on weighted densities}

A differential operator on $\mathbb{R}^{1|n}$ is an operator on
$\mathbb{R}[x,\theta]$ of the form:
\begin{equation}\label{diff}
A=\sum_{k=0}^M
\sum_{\varepsilon=(\varepsilon_1,\cdots,\varepsilon_n)}a_{k,\epsilon}(x,\theta)\partial_x^k
\partial_{1}^{\varepsilon_1}\cdots\partial_{n}^{\varepsilon_n};\,\,
\varepsilon_i=0,1;\,\,M\in\mathbb{N}.
\end{equation}
Of course any differential operator defines a linear mapping
$F\alpha_n^\lambda\mapsto(AF)\alpha_n^\mu$ from
$\mathbb{F}^n_{\lambda}$ to $\mathbb{F}^n_{\mu}$ for any
$\lambda,\,\mu\in\mathbb{R}$, thus the space of differential
operators becomes a family of ${\mathcal K}(n)-$modules $
\mathbb{D}^n_{\lambda,\mu}$ for the natural action:
\begin{equation}\label{d-action}
{X_F}\cdot A=\mathbb{L}^{\mu}_{X_F}\circ A-(-1)^{|A||F|} A\circ
\mathbb{L}^{\lambda}_{X_F}.
\end{equation}

Similarly, we consider  a multi-parameter family of ${
\mathcal{K}}(n)$-modules on the space
$\mathbb{D}^n_{\lambda_1,\dots,\lambda_m;\mu}$ of multi-linear
differential operators: $~A: {\mathbb
F}_{\lambda_1}^n\otimes\cdots\otimes\mathbb{F}_{\lambda_m}^n\longrightarrow{\mathbb
F}_\mu^n$ with the natural $\mathcal{K}(n)-$action:
\begin{equation*}
{X_F}\cdot A={\mathbb L}_{X_F}^\mu\circ A-(-1)^{|A||F|}A\circ
{\mathbb L}_{X_F}^{\lambda_1,\dots,\lambda_m},\end{equation*}where
${\mathbb L}_{X_F}^{\lambda_1,\dots,\lambda_m}$ is defined by the
Leibnitz rule. We also consider the ${\mathcal{K}}(n)-$module
$\Pi\left(\mathbb{D}^n_{\lambda_1,\dots,\lambda_m;\mu}\right)$
with the $\mathcal{K}(n)-$action ($\Pi$ is the change of parity
operator):
\begin{equation*}
{X_F}\cdot \Pi(A)=\Pi\left({\mathbb L}_{X_F}^\mu\circ
A-(-1)^{(|A|+1)|F|}A\circ {\mathbb
L}_{X_F}^{\lambda_1,\dots,\lambda_m}\right).\end{equation*} Since
$-\eta_i^2=\partial_x,$ and $\partial_i=\eta_i-\theta_i\eta_i^2,$
every differential operator $A\in\mathbb{D}^n_{\lambda,\mu}$ can
be expressed in the form
\begin{equation}
\label{diff1}
A(F\alpha_n^\lambda)=\sum_{\ell=(\ell_1,\dots,\ell_n)}a_\ell(x,\theta)
\eta_1^{\ell_1}\dots\eta_n^{\ell_n}(F)\alpha_n^\mu,
\end{equation}
where the coefficients $a_\ell(x,\theta)$ are arbitrary polynomial
functions.
\medskip

The Lie superalgebra ${\mathcal K}(n-1)$ can be realized  as a
subalgebra of ${\mathcal K}(n)$:
\begin{equation*}
 {\mathcal K}(n-1)=\Big\{X_F\in{\mathcal
 K}(n)~|~\partial_{n}F=0\Big\}.
\end{equation*} Therefore,
$\mathbb{D}^n_{\lambda_1,\dots,\lambda_m;\mu}$ and
$\mathbb{F}^n_{\lambda}$ are ${\mathcal K}(n-1)-$modules. Note that, for any  $i$ in
$\{1,2,\dots,n-1\}$, ${\mathcal
K}(n-1)$ is isomorphic to
\begin{equation*}
 {\mathcal K}(n-1)^i=\Big\{X_F\in{\mathcal
 K}(n)~|~\partial_{i}F=0\Big\}.
 \end{equation*}
 \begin{proposition}
\label{iso1} As a $\frak{ aff}(n-1|1)-$module, we have
\begin{align}\label{nizar}
\mathbb{D}^n_{\lambda,\mu;\nu}\simeq\widetilde{\mathbb{D}}^{n-1}_{\lambda,\mu;\nu}:=~~&
{\mathbb{D}}^{n-1}_{\lambda,\mu;\nu}\oplus
{\mathbb{D}}^{n-1}_{\lambda+\frac{1}{2},\mu+\frac{1}{2};\nu}\oplus
{\mathbb{D}}^{n-1}_{\lambda,\mu+\frac{1}{2};\nu+\frac{1}{2}}\oplus
{\mathbb{D}}^{n-1}_{\lambda+\frac{1}{2},\mu;\nu+\frac{1}{2}}\oplus\notag\\[6pt]
&\,
\Pi\left({\mathbb{D}}^{n-1}_{\lambda,\mu;\nu+\frac{1}{2}}\oplus
{\mathbb{D}}^{n-1}_{\lambda,\mu+\frac{1}{2};\nu}\oplus
{\mathbb{D}}^{n-1}_{\lambda+\frac{1}{2},\mu;\nu}\oplus
{\mathbb{D}}^{n-1}_{\lambda+\frac{1}{2},\mu+\frac{1}{2};\nu+\frac{1}{2}}
\right).
\end{align}
\end{proposition}
\begin{proofname}. For any $F\in\mathbb{R}[x,\theta]$, we write $$F=F_1+F_2\theta_n\quad \text{where}\quad\partial_nF_1=\partial_nF_2=0$$ and we prove that
$$\mathbb{L}^\lambda_{X_H}F=\mathbb{L}^\lambda_{X_H}F_1+(\mathbb{L}^{\lambda+\frac{1}{2}}_{X_H}F_2)\theta_n.$$
Thus, it is clear that the map
\begin{equation}\label{varphi}
\begin{array}{lcll} \varphi_\lambda:&\mathbb{F}^n_{\lambda} &\rightarrow&
\mathbb{F}^{n-1}_{\lambda}\oplus\Pi(\mathbb{F}^{n-1}_{\lambda+\frac{1}{2}})\\
&F\alpha_n^\lambda&\mapsto
&(F_1\alpha_{n-1}^\lambda,\,\Pi(F_2\alpha_{n-1}^{\lambda+\frac{1}{2}})),
\end{array}
\end{equation}
is $\frak{ aff}(n-1|1)-$isomorphism.  So, we get the natural
$\frak{ aff}(n-1|1)-$isomorphism  from
$\mathbb{F}^n_{\lambda}\otimes\mathbb{F}^n_{\mu}$ to
\begin{equation*}
 \mathbb{F}^{n-1}_{\lambda}\otimes\mathbb{F}^{n-1}_{\mu}\oplus
\mathbb{F}^{n-1}_{\lambda}\otimes\Pi(\mathbb{F}^{n-1}_{\mu+\frac{1}{2}})
\oplus\Pi(\mathbb{F}^{n-1}_{\lambda+\frac{1}{2}})\otimes\mathbb{F}^{n-1}_{\mu}
\oplus\Pi(\mathbb{F}^{n-1}_{\lambda+\frac{1}{2}})\otimes\Pi(\mathbb{F}^{n-1}_{\mu+\frac{1}{2}})
\end{equation*} denoted $\psi_{\lambda,\mu}$.
Therefore, we deduce a $\frak{ aff}(n-1|1)-$isomorphism:
\begin{equation}\label{Psi}
\begin{array}{llll} \Psi_{\lambda,\mu,\nu}:&\widetilde{\mathbb{D}}^{n-1}_{\lambda,\mu;\nu}
 &\rightarrow&
{\mathbb{D}}^{n}_{\lambda,\mu;\nu}\\[2pt]
&A&\mapsto&\varphi^{-1}_\nu\circ A\circ\psi_{\lambda,\mu}.
\end{array}
\end{equation}Here, we identify
the $\frak{ aff}(n-1|1)-$modules via the following isomorphisms:
\begin{gather*}\begin{array}{llllllll}
\Pi\left({\mathbb{D}}^{n-1}_{\lambda,\mu;\nu'}\right)
&\rightarrow&
\mathrm{Hom_{diff}}\left(\mathbb{F}^{n-1}_\lambda\otimes\mathbb{F}^{n-1}_\mu,
\Pi(\mathbb{F}^{n-1}_{\nu'})\right), &\Pi(A)&\mapsto&\Pi\circ
A,\\[10pt] \Pi\left({\mathbb{D}}^{n-1}_{\lambda,\mu';\nu}\right)
&\rightarrow&
\mathrm{Hom_{diff}}\left(\mathbb{F}^{n-1}_\lambda\otimes
\Pi(\mathbb{F}^{n-1}_{\mu'}), \mathbb{F}^{n-1}_{\nu}\right),
&\Pi(A)&\mapsto& A\circ(1\otimes\Pi),\\[10pt]
\Pi\left({\mathbb{D}}^{n-1}_{\lambda',\mu;\nu}\right)
&\rightarrow&
\mathrm{Hom_{diff}}\left(\Pi(\mathbb{F}^{n-1}_{\lambda'})
\otimes\mathbb{F}^{n-1}_\mu, \mathbb{F}^{n-1}_{\nu}\right),
&\Pi(A)&\mapsto& A\circ(\Pi\otimes\sigma),\\[10pt]
\Pi\left({\mathbb{D}}^{n-1}_{\lambda',\mu';\nu'}\right)
&\rightarrow&
\mathrm{Hom_{diff}}\left(\Pi(\mathbb{F}^{n-1}_{\lambda'})\otimes\Pi(\mathbb{F}^{n-1}_{\mu'}),
\Pi(\mathbb{F}^{n-1}_{\nu'})\right), &\Pi(A)&\mapsto&\Pi\circ
A\circ(\Pi\otimes\sigma\circ\Pi),\\[10pt]
{\mathbb{D}}^{n-1}_{\lambda,\mu';\nu'} &\rightarrow&
\mathrm{Hom_{diff}}\left(\mathbb{F}^{n-1}_{\lambda}\otimes\Pi(\mathbb{F}^{n-1}_{\mu'}),
\Pi(\mathbb{F}^{n-1}_{\nu'})\right),&A&\mapsto&\Pi\circ
A\circ(1\otimes\Pi),\\[10pt]
{\mathbb{D}}^{n-1}_{\lambda',\mu';\nu} &\rightarrow&
\mathrm{Hom_{diff}}\left(\Pi(\mathbb{F}^{n-1}_{\lambda'})\otimes\Pi(\mathbb{F}^{n-1}_{\mu'}),
\mathbb{F}^{n-1}_{\nu}\right),&A&\mapsto&
A\circ(\Pi\otimes\sigma\circ\Pi),\\[10pt]
{\mathbb{D}}^{n-1}_{\lambda',\mu;\nu'} &\rightarrow&
\mathrm{Hom_{diff}}\left(\Pi(\mathbb{F}^{n-1}_{\lambda'})\otimes\mathbb{F}^{n-1}_{\mu},
\Pi(\mathbb{F}^{n-1}_{\nu'})\right),&A&\mapsto&\Pi\circ
A\circ(\Pi\otimes\sigma),
\end{array}
\end{gather*}
where
$\lambda'=\lambda+\frac{1}{2},~\mu'=\mu+\frac{1}{2},~\nu'=\nu+\frac{1}{2}$
and $\sigma(F)=(-1)^{|F|}F$. \end{proofname}\hfill$\Box$

\section{$\frak{aff}(1)-$ and $\frak{aff}(1|1)-$invariant bilinear differential
operators }

In this section we will investigate $\frak{aff}(1)-$ and $\frak{aff}(1|1)-$invariant bilinear
differential operators on tensor densities. These results will be useful for the computation
of cohomology.

\begin{proposition}\label{trans2}\cite{bn}  There exist only the
$\frak{aff}(1)-$invariant bilinear differential operators
\begin{equation*}
J_k^{\tau,\lambda}:
\mathcal{F}_\tau\otimes\mathcal{F}_\lambda\longrightarrow\mathcal{F}_{\tau+\lambda+k},\quad
(\varphi dx^\tau,\phi dx^\lambda)\mapsto
J_k^{\tau,\lambda}(\varphi,\phi)dx^{\tau+\lambda+k}
\end{equation*} given by
\begin{equation*}
J_k^{\tau,\lambda}(\varphi,\phi)=\sum_{0\leq i, j,
i+j=k}c_{i,j}^{\tau,\lambda}\varphi^{(i)}\phi^{(j)},
\end{equation*} where $k\in\mathbb{N}$ and the coefficients $c_{i,j}^{\tau,\lambda}$ are constants.
\end{proposition}

\begin{theorem}\label{main}   There are only the
following $\mathfrak{aff}(1|1)-$invariant bilinear differential
operators acting in the spaces $\mathbb{F}^1_{\lambda}$:
\begin{equation*}\begin{array}{ll}
\mathbb{J}_{k}^{\tau,\lambda;1}:
\mathbb{F}^1_\tau\otimes\mathbb{F}^1_\lambda&\longrightarrow\mathbb{F}^1_{\tau+\lambda+k}\\
(F \alpha^\tau, G \alpha^\lambda)&\mapsto
\mathbb{J}_{k}^{\tau,\lambda;1}(F,G)\alpha^{\tau+\lambda+k},\end{array}
\end{equation*}  where $k\in{1\over2}\mathbb{N}$.
The operators $\mathbb{J}_{k}^{\tau,\lambda;1}$ labeled by
semi-integer $k$ are odd; they are given by
\begin{equation}\label{op11}\begin{array}{lllll}
\mathbb{J}_{k}^{\tau,\lambda;1}(F,G)&=\displaystyle\sum_{i+j=[k]}
\left(\Gamma_{i,j}^{\tau,\lambda;1}{\eta_1}(F^{(i)})G^{(j)}+\widetilde{\Gamma_{i,j}^{\tau,\lambda;1}}
(-1)^{|F|}F^{(i)}{\eta_1}(G^{(j)})\right).
\end{array}\end{equation}
The operators $\mathbb{J}_{k}^{\tau,\lambda;1}$, where $k\in
\mathbb{N}$, are even; set $\mathbb{J}_0^{\tau,\lambda;1}(F,G)=FG$ and
\begin{equation}\label{op12}\begin{array}{lllll}
\mathbb{J}_k^{\tau,\lambda;1}(F,G)&=\displaystyle\sum_{i+j=k}\Upsilon_{i,j}^{\tau,\lambda;1}F^{(i)}G^{(j)}+
\displaystyle\sum_{i+j=k-1}
\widetilde{\Upsilon_{i,j}^{\tau,\lambda;1}}(-1)^{|F|}{\eta_1}(F^{(i)}){\eta_1}(G^{(j)}),\end{array}
\end{equation}where $[k]$ denotes
the integer part of $k$, $k>0$, and
$\Gamma_{i,j}^{\tau,\lambda,1},~\widetilde{\Gamma_{i,j}^{\tau,\lambda,1}},~\Upsilon_{i,j}^{\tau,\lambda,1}$ and $\widetilde{\Upsilon_{i,j}^{\tau,\lambda,1}}$ are constants.
\end{theorem}

\begin{proofname}.  Let $\mathbb{T}^1:
\mathbb{F}^1_\tau\otimes\mathbb{F}^1_\lambda\longrightarrow\mathbb{F}^1_\mu$
be an $\mathfrak{aff}(1|1)-$invariant differential operator. Using
the fact that, as $\mathfrak{vect}(1)-$modules,
\begin{equation}\label{decomposition}\mathbb{F}^1_\tau\otimes\mathbb{F}^1_\lambda\simeq
\mathcal{F}_\tau\otimes\mathcal{F}_\lambda\oplus\Pi(\mathcal{F}_{\tau+\frac{1}{2}}\otimes
\mathcal{F}_{\lambda+\frac{1}{2}})\oplus
\mathcal{F}_\tau\otimes\Pi(\mathcal{F}_{\lambda+\frac{1}{2}})\oplus\Pi(\mathcal{F}_{\tau+\frac{1}{2}})\otimes
\mathcal{F}_\lambda
\end{equation}
and $$
\mathbb{F}^1_\mu\simeq\mathcal{F}_\mu\oplus\Pi(\mathcal{F}_{\mu+\frac{1}{2}}),
$$ we can deduce that the restriction of $\mathbb{T}^1$ to each
component of the right hand side of (\ref{decomposition}) is $\frak{aff}(1)-$invariant.
So, the parameters $\tau,$ $\lambda$ and $\mu$ must
satisfy $$\mu=\lambda+\tau+k,\quad\hbox{ where }\quad
k\in{1\over2}\mathbb{N}.$$ The corresponding operators will be
denoted $\mathbb{J}_{k}^{\tau,\lambda;1}$. Obviously, if $k$ is
integer, then the operator $\mathbb{J}_{k}^{\tau,\lambda;1}$ is even
and its restriction to each component of the right hand side of
(\ref{decomposition}) coincides (up to a scalar factor) with the
respective $\frak{aff}(1)-$invariant operators:

\begin{equation}\label{restric1}
  \begin{cases}
    & \text{${J}_k^{\tau,\lambda}~~~~~~~:
\mathcal{F}_\tau\otimes\mathcal{F}_\lambda\longrightarrow\mathcal{F}_{\mu},$}\\
     & \text{${J}_{k-1}^{\tau+\frac{1}{2},\lambda+\frac{1}{2}}:
\Pi(\mathcal{F}_{\tau+\frac{1}{2}})\otimes\Pi(\mathcal{F}_{\lambda+\frac{1}{2}})\longrightarrow\mathcal{F}_{\mu}$},
\\ & \text{${J}_{k}^{\tau,\lambda+\frac{1}{2}}~~~:
\mathcal{F}_\tau\otimes\Pi(\mathcal{F}_{\lambda+\frac{1}{2}})\longrightarrow\Pi(\mathcal{F}_{\mu+\frac{1}{2}}$}),
\\ & \text{${J}_{k}^{\tau+\frac{1}{2},\lambda}~~~:
\Pi(\mathcal{F}_{\tau+\frac{1}{2}})\otimes\mathcal{F}_\lambda\longrightarrow\Pi(\mathcal{F}_{\mu+\frac{1}{2}}$}).
\\
  \end{cases}
  \end{equation}
If $k$ is semi-integer, then the operator
$\mathbb{J}_{k}^{\tau,\lambda;1}$ is odd and its  restriction to each
component of the right hand side of (\ref{decomposition})
coincides (up to a scalar factor ) with the respective
$\frak{aff}(1)-$invariant operators:
\begin{equation}\label{restric2}
  \begin{cases}
    & \text{$J_{[k]+1}^{\tau,\lambda}~~~~~~~:
\mathcal{F}_\tau\otimes\mathcal{F}_\lambda\longrightarrow\Pi(\mathcal{F}_{\mu+\frac{1}{2}}),$}\\
     & \text{$J_{[k]}^{\tau+\frac{1}{2},\lambda+\frac{1}{2}}:
\Pi(\mathcal{F}_{\tau+\frac{1}{2}})\otimes\Pi(\mathcal{F}_{\lambda+\frac{1}{2}})
\longrightarrow\Pi(\mathcal{F}_{\mu+\frac{1}{2}})$}, \\ &
\text{$J_{[k]}^{\tau,\lambda+\frac{1}{2}}~~~:
\mathcal{F}_\tau\otimes\Pi(\mathcal{F}_{\lambda+\frac{1}{2}})\longrightarrow\mathcal{F}_{\mu}$},
\\ & \text{$J_{[k]}^{\tau+\frac{1}{2},\lambda}~~~:
\Pi(\mathcal{F}_{\tau+\frac{1}{2}})\otimes\mathcal{F}_\lambda\longrightarrow\mathcal{F}_{\mu}$}.
\\
  \end{cases}
  \end{equation}
More precisely, let $F\alpha^\tau\otimes G
\alpha^\lambda\in\mathbb{F}^1_\tau\otimes\mathbb{F}^1_\lambda$, where
$F=f_0+\theta_1 f_1$ and $G=g_0+\theta_1 g_1$, with
$f_0,\,f_1,\,g_0,\,g_1\in\mathbb{R}[x]$. Then if $k$ is integer, we
have
\begin{equation}\label{integer}
\mathbb{J}_{k}^{\tau,\lambda;1}(\varphi,\psi)=\Big[a_1J_k^{\tau,\lambda}(f_0,g_0)+
a_2J_{k-1}^{\tau+\frac{1}{2},\lambda+\frac{1}{2}}(f_1,g_1)
+\theta_1\left(a_3J_{k}^{\tau,\lambda+\frac{1}{2}}(f_0,g_1)+a_4J_{k}^{\tau+\frac{1}{2},\lambda}(f_1,g_0)
\right)\Big]\alpha_1^\mu
\end{equation}
and if $k$ is semi-integer, we have
\begin{equation}\label{semi}
\mathbb{J}_{k}^{\tau,\lambda;1}(\varphi,\psi)=\Big[b_1J_{[k]}^{\tau,\lambda+\frac{1}{2}}(f_0,g_1)+
b_2J_{[k]}^{\tau+\frac{1}{2},\lambda}(f_1,g_0)
+\theta_1\left(b_3J_{[k]+1}^{\tau,\lambda}(f_0,g_0)+b_4J_{[k]}^{\tau+\frac{1}{2},\lambda+\frac{1}{2}}(f_1,g_1)
\right)\Big]\alpha_1^\mu,
\end{equation}where the $a_i$ and $b_i$ are constants.
The invariance of $\mathbb{J}_{k}^{\tau,\lambda;1}$ with respect to
$X_{\theta_1}$ reads:
\begin{equation}\label{T1}
{\frak L}_{X_{\theta_1}}^\mu\circ
\mathbb{J}_{k}^{\tau,\lambda;1}-(-1)^{2k}\mathbb{J}_{k}^{\tau,\lambda;1}\circ
{\frak L}_{X_{\theta_1}}^{(\tau,\lambda)}=0.\end{equation}

The formula (\ref{T1}) allows us to determine the coefficients
$a_i$ and $b_i$. More precisely, the invariance property with
respect to $X_{\theta_1}$ yields:
\begin{enumerate}
\item If $k$ is integer
\begin{equation*}
\begin{array}{ccc}
c_{i,k-i}^{\tau,\lambda+\frac{1}{2}}&=&c_{i-1,k-i}^{\tau+\frac{1}{2},\lambda+\frac{1}{2}}+
c_{i,k-i}^{\tau,\lambda}\\[8pt]
c_{i,k-i}^{\tau+\frac{1}{2},\lambda}&=&c_{i,k-i}^{\tau,\lambda}-c_{i,k-i-1}^{\tau+\frac{1}{2},\lambda+\frac{1}{2}}
\end{array}
\end{equation*}
\item If $k$ is semi-integer
\begin{equation*}
\begin{array}{ccc}
c_{i,k-i+\frac{1}{2}}^{\tau,\lambda}&=&-c_{i,k-i-\frac{1}{2}}^{\tau,\lambda+\frac{1}{2}}-
c_{i-1,k-i+\frac{1}{2}}^{\tau+\frac{1}{2},\lambda}\\[8pt]
c_{i+1,k-i-\frac{1}{2}}^{\tau,\lambda}&=&-c_{i+1,k-i-\frac{3}{2}}^{\tau,\lambda+\frac{1}{2}}-
c_{i,k-i-\frac{1}{2}}^{\tau+\frac{1}{2},\lambda}\\[8pt]
c_{i,k-i-\frac{1}{2}}^{\tau+\frac{1}{2},\lambda+\frac{1}{2}}&=&-c_{i,k-i-\frac{1}{2}}^{\tau,\lambda+\frac{1}{2}}+
c_{i,k-i-\frac{1}{2}}^{\tau+\frac{1}{2},\lambda},
\end{array}
\end{equation*}
\end{enumerate}
where $c_{i,j}^{\tau,\lambda}=0~\forall i<0~\hbox{or}~j<0.$
Therefore, we easily check that $\mathbb{J}_{k}^{\tau,\lambda;1}$ is expressed as in
(\ref{op11}-\ref{op12}).
\end{proofname}

\section{The $\frak{aff}(2|1)-$invariant bilinear differential operators }

In this section we will investigate $\frak{aff}(2|1)-$invariant bilinear
differential operators on tensor densities. These results allow us to compute
$\frak{aff}(2|1)-$relative cohomology.

\begin{theorem}\label{main}   The space of $\mathfrak{aff}(2|1)-$invariant bilinear differential
operators acting in the spaces $\mathbb{F}^2_{\lambda}$:
\begin{equation*}\begin{array}{ll}
\mathbb{J}_{k}^{\tau,\lambda;2}:
\mathbb{F}^2_\tau\otimes\mathbb{F}^2_\lambda&\longrightarrow\mathbb{F}^2_{\tau+\lambda+k}\\
(F \alpha^\tau, G \alpha^\lambda)&\mapsto
\mathbb{J}_{k}^{\tau,\lambda;2}(F,G)\alpha^{\tau+\lambda+k},\end{array}
\end{equation*}  where $k\in{1\over2}\mathbb{N}$, is purely even and it is spanned by the operator $\mathbb{J}_0^{\tau,\lambda;2}(F,G)=FG$ for $k=0,$ the operator
\begin{equation}\label{opp23}\begin{array}{lllll}
\mathbb{J}_1^{\tau,\lambda;2}(F,G)&=\displaystyle\sum_{i+j=1}\Gamma_{i,j,1}^{\tau,\lambda;2}F^{(i)}G^{(j)}+
\Gamma_{0,0,2}^{\tau,\lambda;2}(-1)^{|F|}({\eta_1}(F){\eta_1}(G)+
{\eta_2}(F){\eta_2}(G))+\\&
\Gamma_{0,0,3}^{\tau,\lambda;2}(-1)^{|F|}({\eta_1}(F){\eta_2}(G)-{\eta_2}(F){\eta_1}(G))+
\Gamma_{0,0,4}^{\tau,\lambda;2}{\eta_1}{\eta_2}(F)G+
\Gamma_{0,0,5}^{\tau,\lambda;2}F{\eta_1}{\eta_2}(G)\end{array}
\end{equation}
for $k=1,$
and the operator
\begin{equation}\label{opp24}\begin{array}{lllll}
\mathbb{J}_k^{\tau,\lambda;2}(F,G)&=\displaystyle\sum_{i+j=k}\Gamma_{i,j,1}^{\tau,\lambda;2}F^{(i)}G^{(j)}+
\displaystyle\sum_{i+j=k-1}\Gamma_{i,j,2}^{\tau,\lambda;2}(-1)^{|F|}({\eta_1}(F^{(i)}){\eta_1}(G^{(j)})+
{\eta_2}(F^{(i)}){\eta_2}(G^{(j)}))\\&+
\displaystyle\sum_{i+j=k-1}\Gamma_{i,j,3}^{\tau,\lambda;2}(-1)^{|F|}({\eta_1}(F^{(i)}){\eta_2}(G^{(j)})-
{\eta_2}(F^{(i)}){\eta_1}(G^{(j)}))+
\displaystyle\sum_{i+j=k-1}\Gamma_{i,j,4}^{\tau,\lambda;2}{\eta_1}{\eta_2}(F^{(i)})G^{(j)}\\&+
\displaystyle\sum_{i+j=k-1}\Gamma_{i,j,5}^{\tau,\lambda;2}F^{(i)}{\eta_1}{\eta_2}(G^{(j)})+
\displaystyle\sum_{i+j=k-2}\Gamma_{i,j,6}^{\tau,\lambda;2}{\eta_1}{\eta_2}(F^{(i)}){\eta_1}{\eta_2}(G^{(j)})
\end{array}
\end{equation}

for $k\geq2;$ where $\Gamma_{i,j,s}^{\tau,\lambda,2},~s\in \{1,\ldots,6\},$ are constants.
\end{theorem}

\begin{proofname}.
Let $\mathbb{T}^2:
\mathbb{F}^{2}_\tau\otimes\mathbb{F}^{2}_\lambda\longrightarrow\mathbb{F}^{2}_\mu$
be an $\mathfrak{aff}(2|1)-$invariant bilinear differential
operator. Observe that the $\mathfrak{aff}(2|1)-$invariance of
$\mathbb{T}^2$ is equivalent to invariance with respect to $X_{\theta_1\theta_2}$ and the
subsuperalgebras $\mathfrak{aff}(1|1)$ and $\mathfrak{aff}(1|1)_1.$
Note that, the $\mathfrak{aff}(1|1)-$invariant elements of
$\Pi(\mathbb{D}^{2}_{\tau,\lambda;\mu})$ can be deduced from those
given in (\ref{op11}) and (\ref{op12}) by using  the
following $\mathfrak{aff}(1|1)-$isomorphism
\begin{equation}
\label{od} \mathbb{D}^{1}_{\tau,\lambda;\mu}\rightarrow
\Pi(\mathbb{D}^{1}_{\tau,\lambda;\mu}),\quad~A\mapsto
\Pi(A\circ(\sigma\otimes\sigma)).
\end{equation}
Now, by isomorphism (\ref{Psi}) we exhibit the
$\mathfrak{aff}(1|1)-$invariant elements of
$\mathbb{D}^{2}_{\tau,\lambda;\mu}$. Of course, these elements are
identically zero if $2(\mu-\tau-\lambda)\notin\mathbb{N}$. So, the
parameters $\tau,$ $\lambda$ and $\mu$ must satisfy
\[
\mu=\tau+\lambda+k,\quad\hbox{ where }\quad
k\in{1\over2}\mathbb{N}.
\]
The corresponding operators will be denoted
$\mathbb{J}_{k}^{\tau,\lambda,2}$. Obviously, if $k$ is integer,
then the operator $\mathbb{J}_{k}^{\tau,\lambda,2}$ is even and if
$k$ is semi-integer, then the operator
$\mathbb{J}_{k}^{\tau,\lambda,2}$ is odd.

Any $\mathfrak{aff}(2|1)-$invariant element
$\mathbb{J}_{k}^{\tau,\lambda,2}$ of $\mathbb{D}^{2}_{\tau,
\lambda;\mu}$ with $\mu=\tau+\lambda+k$, can be expressed as
follows: {\small\begin{equation} \label{nadhir}
\begin{array}{lllll}
\mathbb{J}_{k}^{\tau,\lambda,2}=
&\Psi_{\tau,\lambda,\mu}
\left(\mathbb{J}_{k}^{\tau,\lambda,1}+
\mathbb{J}_{k-1}^{\tau+{1\over2},\lambda+{1\over2},1}+
\mathbb{J}_{k}^{\tau,\lambda+{1\over2},1}+
\mathbb{J}_{k}^{\tau+{1\over2},\lambda,1}
\right)+\\[5pt]
 &\Psi_{\tau,\lambda,\mu}
\left(\Pi\left(\left(\mathbb{J}_{k+{1\over2}}^{\tau,\lambda,1}+
\mathbb{J}_{k-{1\over2}}^{\tau,\lambda+{1\over2},1}+
\mathbb{J}_{k-{1\over2}}^{\tau+{1\over2},\lambda,1}+
\mathbb{J}_{k-{1\over2}}^{\tau+{1\over2},\lambda+{1\over2},1}\right)\circ
(\sigma\otimes\sigma)\right)
\right)
\end{array}
\end{equation}} where $\mathbb{J}_{k}^{\tau,\lambda,1}$ are defined
by (\ref{op11}--\ref{op12}). The
invariance of $\mathbb{J}_{k}^{\tau,\lambda,2}$ with respect to
$\mathfrak{aff}(1|1)_1$ imposes some supplementary conditions over
the coefficients of the operators $\mathbb{J}_{k}^{\tau,\lambda,2}$. By a direct computation, we
get:

$\bullet$ For $k\in\mathbb{N}+\frac{1}{2},$
{\small\[\begin{array}{llllll}
\widetilde{\Gamma_{i,k-i-\frac{1}{2}}^{\tau,\lambda+\frac{1}{2};1}}&=&\widetilde{\Gamma_{i,k-i-\frac{1}{2}}^{\tau,\lambda;1}}
+\widetilde{\Gamma_{i-1,k-i-\frac{1}{2}}^{\tau+\frac{1}{2},\lambda+\frac{1}{2};1}},\quad
&\Upsilon_{i,k-i-\frac{1}{2}}^{\tau+\frac{1}{2},\lambda+\frac{1}{2};1}&=&-\Upsilon_{i,k-i-\frac{1}{2}}^{\tau,\lambda+\frac{1}{2};1}
+\Upsilon_{i,k-i-\frac{1}{2}}^{\tau+\frac{1}{2},\lambda;1},\\[10pt]
\Gamma_{i,k-i-\frac{1}{2}}^{\tau,\lambda+\frac{1}{2};1}&=&\Gamma_{i,k-i-\frac{1}{2}}^{\tau,\lambda;1}
-\Gamma_{i-1,k-i-\frac{1}{2}}^{\tau+\frac{1}{2},\lambda+\frac{1}{2};1},\quad
&\widetilde{\Upsilon_{i,k-i-\frac{3}{2}}^{\tau+\frac{1}{2},\lambda+\frac{1}{2};1}}&=&-\widetilde{\Upsilon_{i,k-i-\frac{3}{2}}^{\tau,\lambda+\frac{1}{2};1}}
+\Upsilon_{i,k-i-\frac{3}{2}}^{\tau+\frac{1}{2},\lambda;1},\\[10pt]
\widetilde{\Gamma_{i,k-i-\frac{1}{2}}^{\tau+\frac{1}{2},\lambda;1}}&=&-\widetilde{\Gamma_{i,k-i-\frac{1}{2}}^{\tau,\lambda;1}}
+\widetilde{\Gamma_{i,k-i-\frac{3}{2}}^{\tau+\frac{1}{2},\lambda+\frac{1}{2};1}},\quad
&\Upsilon_{i,k-i+\frac{1}{2}}^{\tau,\lambda;1}&=&\Upsilon_{i,k-i-\frac{1}{2}}^{\tau,\lambda+\frac{1}{2};1}
+\Upsilon_{i-1,k-i+\frac{1}{2}}^{\tau+\frac{1}{2},\lambda;1},\\[10pt]
\Gamma_{i,k-i-\frac{1}{2}}^{\tau+\frac{1}{2},\lambda;1}&=&\Gamma_{i,k-i-\frac{1}{2}}^{\tau,\lambda;1}
+\Gamma_{i,k-i-\frac{3}{2}}^{\tau+\frac{1}{2},\lambda+\frac{1}{2};1},\quad
&\widetilde{\Upsilon_{i,k-i-\frac{1}{2}}^{\tau,\lambda;1}}&=&
-\widetilde{\Upsilon_{i,k-i-\frac{3}{2}}^{\tau,\lambda+\frac{1}{2};1}}
+\widetilde{\Upsilon_{i-1,k-i-\frac{1}{2}}^{\tau+\frac{1}{2},\lambda;1}},
\end{array}
\]}
where $\Gamma_{i,j}^{\tau,\lambda}=\widetilde{\Gamma_{i,j}^{\tau,\lambda}}=\Upsilon_{i,j}^{\tau,\lambda}=
\widetilde{\Upsilon_{i,j}^{\tau,\lambda}}=0~\forall i<0~\hbox{or}~j<0.$

$\bullet$ For $k\in\mathbb{N},$
{\small\[\begin{array}{llllll}
\Upsilon_{i,k-i}^{\tau,\lambda+\frac{1}{2};1}&=&\Upsilon_{i,k-i}^{\tau,\lambda;1}
-\widetilde{\Upsilon_{i-1,k-i}^{\tau,\lambda;1}},\quad
&\widetilde{\Gamma_{i,k-i}^{\tau,\lambda;1}}&=&\widetilde{\Gamma_{i,k-i-1}^{\tau,\lambda+\frac{1}{2};1}}
-\widetilde{\Gamma_{i-1,k-i}^{\tau+\frac{1}{2},\lambda;1}},\\[10pt]
\Upsilon_{i,k-i}^{\tau+\frac{1}{2},\lambda+\frac{1}{2};1}&=&\Upsilon_{i,k-i}^{\tau,\lambda;1}
+\widetilde{\Upsilon_{i,k-i-1}^{\tau,\lambda;1}},\quad
&\widetilde{\Gamma_{i,k-i}^{\tau,\lambda;1}}&=&\widetilde{\Gamma_{i,k-i-1}^{\tau,\lambda+\frac{1}{2};1}}
+\Gamma_{i-1,k-i}^{\tau+\frac{1}{2},\lambda;1},\\[10pt]
\widetilde{\Upsilon_{i,k-i-1}^{\tau,\lambda+\frac{1}{2};1}}&=&-\Upsilon_{i,k-i-1}^{\tau+\frac{1}{2},\lambda+\frac{1}{2};1}
-\widetilde{\Upsilon_{i-1,k-i-1}^{\tau+\frac{1}{2},\lambda+\frac{1}{2};1}},\quad
&\widetilde{\Gamma_{i,k-i-1}^{\tau+\frac{1}{2},\lambda+\frac{1}{2};1}}&=&
\widetilde{\Gamma_{i,k-i-1}^{\tau,\lambda+\frac{1}{2};1}}
-\widetilde{\Gamma_{i,k-i-1}^{\tau+\frac{1}{2},\lambda;1}},\\[10pt]
\widetilde{\Upsilon_{i,k-i-1}^{\tau+\frac{1}{2},\lambda;1}}&=&-\Upsilon_{i,k-i-1}^{\tau+\frac{1}{2}
,\lambda+\frac{1}{2};1}
+\widetilde{\Upsilon_{i,k-i-2}^{\tau+\frac{1}{2},\lambda+\frac{1}{2};1}},\quad
&\widetilde{\Gamma_{i,k-i-1}^{\tau+\frac{1}{2},\lambda+\frac{1}{2};1}}&=&
-\widetilde{\Gamma_{i,k-i-1}^{\tau,\lambda+\frac{1}{2};1}}
+\widetilde{\Gamma_{i,k-i-1}^{\tau+\frac{1}{2},\lambda;1}},
\end{array}
\]}
where $\Gamma_{i,j}^{\tau,\lambda}=\widetilde{\Gamma_{i,j}^{\tau,\lambda}}=\Upsilon_{i,j}^{\tau,\lambda}=
\widetilde{\Upsilon_{i,j}^{\tau,\lambda}}=0~\forall i<0~\hbox{or}~j<0.$
\vskip0.2cm
Finally, the invariance with respect to $X_{\theta_1\theta_2}$ completely  determines
the space of $\mathfrak{aff}(2|1)-$invariant elements of
$\mathbb{D}^{2}_{\tau,\lambda;\mu}.$
\end{proofname}

\section{The $\mathfrak{aff}(n|1)-$Relative Cohomology of $\mathcal{K}(n)$
Acting on $\mathbb{D}^n_{\lambda,\mu}$}

Let us first recall some fundamental concepts from cohomology
theory~(see, e.g., \cite{NIS,C,Fu}).
\subsection{Lie superalgebra cohomology}

 Let $\frak{g}=\frak{g}_{\bar
0}\oplus \frak{g}_{\bar 1}$ be a Lie superalgebra acting on a
superspace $V=V_{\bar 0}\oplus V_{\bar 1}$ and let $\mathfrak{h}$
be a subalgebra of $\mathfrak{g}$. (If $\frak{h}$ is omitted it is
assumed to be $\{0\}$.) The space of $\frak h$-relative
$n$-cochains of $\frak{g}$ with values in $V$ is the
$\frak{g}$-module
\begin{equation*}
C^n(\frak{g},\frak{h}; V ) := \mathrm{Hom}_{\frak
h}(\Lambda^n(\frak{g}/\frak{h});V).
\end{equation*}
The {\it coboundary operator} $ \delta_n: C^n(\frak{g},\frak{h}; V
)\longrightarrow C^{n+1}(\frak{g},\frak{h}; V )$ is a
$\frak{g}$-map satisfying $\delta_n\circ\delta_{n-1}=0$. The
kernel of $\delta_n$, denoted $Z^n(\mathfrak{g},\frak{h};V)$, is
the space of $\frak{h}$-relative $n$-{\it cocycles}, among them,
the elements in the range of $\delta_{n-1}$ are called $\frak{
h}$-relative $n$-{\it coboundaries}. We denote
$B^n(\mathfrak{g},\frak{h};V)$ the space of $n$-coboundaries.

By definition, the $n^{th}$ $\frak{h}$-relative  cohomolgy space is
the quotient space
\begin{equation*}
\mathrm{H}^n
(\mathfrak{g},\frak{h};V)=Z^n(\mathfrak{g},\frak{h};V)/B^n(\mathfrak{g},\frak{h};V).
\end{equation*}
We will only need the formula of $\delta_n$ (which will be simply
denoted $\delta$) in degrees 0 and 1: for $v \in
C^0(\frak{g},\,\frak{h}; V) =V^{\frak h}$,~ $\delta v(g) : =
(-1)^{|g||v|}g\cdot v$, where
\begin{equation*}
V^{\frak h}=\{v\in V~\mid~h\cdot v=0\quad\text{ for all }
h\in\frak h\},
\end{equation*}
and  for  $ \Upsilon\in C^1(\frak{g}, \frak{h};V )$,
\begin{equation*}\delta(\Upsilon)(g,\,h):=
(-1)^{|g||\Upsilon|}g\cdot
\Upsilon(h)-(-1)^{|h|(|g|+|\Upsilon|)}h\cdot
\Upsilon(g)-\Upsilon([g,~h])\quad\text{for any}\quad g,h\in
\frak{g}.
\end{equation*}

\subsection{The space
$\mathrm{H}^1_{\mathrm{diff}}(\mathcal{K}(n),\mathfrak{aff}(n|1);\mathbb{D}^n_{\lambda,\mu})$}

In this  subsection, we will compute the first differential
cohomology spaces  ${\mathrm H}^1_{\rm
diff}(\mathcal{K}(n),\mathfrak{aff}(n|1);
\mathbb{D}^n_{\lambda,\mu})$ for $n=0,~1$ and $2.$ Our main tools are the following two results

\begin{lemma}\label{aff} Any $1$-cocycle $\Upsilon\in
Z^1_{\mathrm{diff}}(\mathcal{K}(n);\mathbb{D}^n_{\lambda,\mu})$
vanishing on $\frak{aff}(n|1)$ is $\mathfrak{aff}(n|1)-$invariant.
\end{lemma}
\begin{proofname}. The $1$-cocycle relation of $\Upsilon$ reads:
\begin{equation}\label{osp1}
(-1)^{|F||\Upsilon|}\mathbb{L}_{X_F}^{\lambda,\mu}
\Upsilon(X_G)-(-1)^{|G|(|F|+|\Upsilon|)}\mathbb{L}_{X_G}^{\lambda,\mu}
\Upsilon(X_F)-\Upsilon([X_F,~X_G])=0,
\end{equation}
where $X_F,\,X_G\in ~\mathcal{K}(n).$  Thus, if $\Upsilon(X_F)=0$
for all $X_F\in\frak{aff}(n|1)$, the equation (\ref{osp1}) becomes
\begin{equation}\label{osp2}
(-1)^{|F||\Upsilon|}\mathbb{L}_{X_F}^{\lambda,\mu}
\Upsilon(X_G)-\Upsilon([X_F,~X_G])=0
\end{equation}
 expressing the $\frak{aff}(n|1)-$invariance  of
$\Upsilon$. 
\end{proofname}

\medskip

\begin{thm}\cite{BIIO}\label{LNI}
Let $ \mathcal{A}^n_{\lambda,
\mu}:\mathfrak{F}^n_\lambda\rightarrow\mathfrak{F}^n_\mu,
\,(F\alpha_{n}^\lambda)\mapsto \mathcal{A}^n_{\lambda,
\mu}(F)\alpha_{n}^{\mu}$ be a non-zero
$\mathfrak{aff}(n|1)-$invariant linear differential operator. Then,
up to a scalar factor, the map $\mathcal{A}^n_{\lambda,\mu}$ is
given by:

\begin{equation} \label{k1}
\begin{array}{llllllllllll}
\mathcal{A}^n_{\lambda,\lambda+k}(F)&=&
F^{(k)},&\hbox{for}~k\in\mathbb{N}~\hbox{and}~~n\in\mathbb{N}\\[5pt]
\mathcal{A}^n_{\lambda-\frac{n}{2},\lambda+k}(F)&=&
\overline{\eta}_{1}\overline{\eta}_{2}\cdots\overline{\eta}_{n}(F^{(k)}),
&\hbox{for}~k\in\mathbb{N}~\hbox{and}~~n~\geq1.
\end{array}
\end{equation}
\end{thm}

\subsubsection{The space $\mathrm{H}^1_{\rm diff}(\mathfrak{vect}(1),\mathfrak{aff}(1);
\mathrm{D}_{\lambda,\mu})$}

The main result of this subsection is the following

\begin{theorem}\label{FirstSect}
\begin{equation}
\label{CohoSpace2} \mathrm{H}^1_{\rm diff}(\mathfrak{vect}(1),\mathfrak{aff}(1);
\mathrm{D}_{\lambda,\mu})\simeq\left\{
\begin{array}{ll}
\mathbb{R} & \hbox{if} ~~\left\{
\begin{array}{l}
 \mu-\lambda=2,3,4 \hbox{ for all
}\lambda,\\[2pt]\lambda=0\hbox{ ~and~
}\mu=1 ,\\[2pt] \lambda=0  \hbox{ or }
\lambda=-4\hbox{ ~and~ }\mu-\lambda=5, \\[2pt]
\lambda=-\frac{5\pm \sqrt{19}}{2}\hbox{ ~and~ }\mu-\lambda
=6,\end{array}
\right.\\[16pt] 0 &\hbox{ otherwise. }
\end{array}
\right.
\end{equation}
For $X\frac{d}{dx}\in\mathfrak{vect}(1)$ and
$f{dx}^{\lambda}\in{\cal F}_\lambda$, we write
\begin{align*}\begin{array}{llll}
C_{\lambda,\lambda+k
}(X\frac{d}{dx})(f{dx}^{\lambda})=C_{\lambda,\lambda+k
}(X,f){dx}^{\lambda+k}. \end{array}\end{align*} The spaces
${\mathrm H}^1_{\rm dif\/f}(\mathfrak{vect}(1),\mathfrak{aff}(1);
\mathrm{D}_{\lambda,\lambda+k})$ are generated by the cohomology
classes of the following 1-cocycles:
\begin{equation}\label{cocycles}\begin{array}{llllllllll}
C_{0,1}(X,f)&=&X''f \\
   C_{\lambda,\lambda+2}(X,f)&=&X^{(3)}f+2X''f' \\
  C_{\lambda,\lambda+3}(X,f)&=&X^{(3)}f'+X''f'' \\
  C_{\lambda,\lambda+4}(X,f)&=&-\lambda
X^{(5)}f+X^{(4)}f'-6X^{(3)}f''-4X''f^{(3)}\\
C_{0,5}(X,f)&=&2X^{(5)}f'-5X^{(4)}f''+10X^{(3)}f^{(3)}+5X''f^{(4)}\\
C_{-4,1}(X,f)&=&12X^{(6)}f+22X^{(5)}f'+5X^{(4)}f''-10X^{(3)}f^{(3)}-5X''f^{(4)}\\
  C_{a_i,a_i+6}(X,f)&=&\alpha_i X^{(7)}f-\beta_i X^{(6)}f'-\gamma_i
X^{(5)}f''-
  5X^{(4)}f^{(3)}+5X^{(3)}f^{(4)}+~&2X''f^{(5)},
\end{array}
\end{equation}
where\begin{equation*}
\begin{array}{llllllll}a_1=-\frac{5+ \sqrt{19}}{2},
&\alpha_1=-\frac{22+ 5\sqrt{19}}{4}, &\beta_1=\frac{31+
7\sqrt{19}}{2}, &\gamma_1=\frac{25+ 7\sqrt{19}}{2}\\[4pt]
a_2=-\frac{5- \sqrt{19}}{2}, &\alpha_2=-\frac{22- 5\sqrt{19}}{4},
& \beta_2=\frac{31- 7\sqrt{19}}{2},&\gamma_2=\frac{25-
7\sqrt{19}}{2}.\end{array}\end{equation*}
\end{theorem}

\begin{proofname}.
Note that, by Lemma \ref{aff}, the $\mathfrak{aff}(1)-$relative cocycles are $\mathfrak{aff}(1)-$invariant bilinear
differential operators. On the other hand, Feigin and Fuchs \cite{ff} calculated
$\mathrm{H}^1_{\rm diff}(\mathfrak{vect}(1);
\mathrm{D}_{\lambda,\mu})$. The result is as follows
\begin{equation}
\label{CohSpace2} \mathrm{H}^1_{\rm diff}(\mathfrak{vect}(1);
\mathrm{D}_{\lambda,\mu})\simeq\left\{
\begin{array}{ll}
\mathbb{R}&\hbox{ if }~~ \mu-\lambda=0,2,3,4 \hbox{ for all
}\lambda,\\[2pt] \mathbb{R}^2& \hbox{ if }~~\lambda=0\hbox{ ~and~
}\mu=1 ,\\[2pt] \mathbb{R}&\hbox{ if }~~ \lambda=0  \hbox{ or }
\lambda=-4\hbox{ ~and~ }\mu-\lambda=5, \\[2pt] \mathbb{R}& \hbox{
if }~~ \lambda=-\frac{5\pm \sqrt{19}}{2}\hbox{ ~and~ }\mu-\lambda
=6,\\[2pt] 0 &\hbox{ otherwise. }
\end{array}
\right.
\end{equation}
Hence, the cohomology spaces $\mathrm{H}^1_{\rm
diff}(\mathfrak{vect}(1),\mathfrak{aff}(1);\mathrm{D}_{\lambda,\mu})$ vanish for $\mu-\lambda\geq 7.$
So we only have to study the $\mathfrak{aff}(1)-$invariant bilinear differential operators $J_k^{-1,\lambda}$ for $k\leq7.$
More precisely, the $1$-cocycle condition imposes conditions on the constants $c_{i,j}:$ we get a linear system
for $c_{i,j}.$ Thereafter, taking into account these conditions, we eliminate all constants underlying
coboundaries. A straightforward but long computation leads to the result.
\end{proofname}

\subsubsection{The space $\mathrm{H}^1_{\rm
diff}(\mathcal{K}(1),\mathfrak{aff}(1|1);\mathbb{D}^1_{\lambda,\mu})$}

The main result of this subsection is the following

\begin{theorem}
\label{th1} \begin{equation}
\label{CohoSpace3}\mathrm{H}^1_{\rm
diff}(\mathcal{K}(1),\mathfrak{aff}(1|1);\mathbb{D}^1_{\lambda,\mu}) \simeq\left\{
\begin{array}{ll}
\mathbb{R} & \hbox{if} ~~\left\{
\begin{array}{llll}
\mu-\lambda=\frac{1}{2}&\hbox{ for  }~\lambda=0,\\[4pt]
\mu-\lambda=\frac{3}{2},2,\frac{5}{2}&\hbox{ for all }~\lambda,\\[4pt] \mu-\lambda=3 &\hbox{ for }
\lambda=0,\,-\frac{5}{2},\\[4pt] \mu-\lambda=4 &\hbox{ for }
\lambda=\frac{-7\pm\sqrt{33}}{4}.
\end{array}
\right.\\[16pt] 0 &\hbox{ otherwise. }
\end{array}
\right.
\end{equation}
\vskip0.2cm
The space $\mathrm{H}^1_{\rm
diff}(\mathcal{K}(1),\mathfrak{aff}(1|1);\mathbb{D}^1_{\lambda,\mu})$ is spanned by the cohomology classes of the
following $1$-cocycles:
\begin{equation*}\label{maincocyc}
  \begin{array}{lllllllllll}
\Upsilon^{1}_{0,\frac{1}{2}}(X_G)(F)&=&
{\eta_1}(G')F\alpha_1^{{1\over2}}.\\[4pt]
\Upsilon^{1}_{\lambda,\lambda+\frac{3}{2}}(X_G)(F\alpha_1^{\lambda})&=&\left\{\begin{array}{ll}{\eta_1}(G'')F
\alpha_1^{\lambda+\frac{3}{2}}
&\hfill\text{
if } \lambda\neq-\frac{1}{2},\\[4pt]\Big({\eta_1}(G'')F
+{\eta_1}(G')F'+(-1)^{|G|}G''{\eta_1}(F)\Big)\alpha_1&\hfill\text{
if }\lambda=-{1\over2}.
 \end{array}\right.\\[8pt]
 \Upsilon^{1}_{\lambda,\lambda+\frac{5}{2}}(X_G)(F\alpha_1^{\lambda})&=&\left\{\begin{array}{ll}\Big( 3{\eta_1}(G'')F'-2\lambda{\eta_1}(G''')F+(-1)^{|G|}(G'''{\eta_1}(F)
\Big)\alpha_1^{\lambda+\frac{5}{2}}
&\hfill\text{
if } \lambda\neq-1,\\[4pt]\Big((-1)^{|G|}(G'''{\eta_1}(F)+2G''{\eta_1}(F'))+
2{\eta_1}(G'')F'+{\eta_1}(G')F''\Big)\alpha_1^{\frac{3}{2}}&\hfill\text{
if }\lambda=-1.
 \end{array}\right.
\\[8pt]
\Upsilon^{1}_{\lambda,\lambda+2}(X_G)(F\alpha_1^{\lambda})&=&
\Big( \frac{2}{3}\lambda G'''F-(-1)^{|G|}{\eta_1}(G''){\eta_1}(F)\Big)\alpha_1^{\lambda+2}\hfill\text{if
}\lambda\neq-1.\\[8pt]
\Upsilon^{1}_{\lambda,\lambda+3}(X_G)(F\alpha_1^{\lambda})&=&
\Big( (-1)^{|G|}{\eta_1}(G''){\eta_1}(F')-\frac{2\lambda+1}{3}\big((-1)^{|G|}{\eta_1}(G'''){\eta_1}(F)
+G'''F'\big)+\\[4pt]&&\frac{\lambda(2\lambda+1)}{6}G^{(4)}F
\Big)\alpha_1^{\lambda+3}\hfill\text{ if } \lambda=0,\frac{-5}{2}.\\[8pt]
\Upsilon^{1}_{\lambda,\lambda+4}(X_G)(F\alpha_1^{\lambda})&=&\Big( (-1)^{|G|}{\eta_1}(G''){\eta_1}(F'')-\frac{2(\lambda+1)}{3}\big(2(-1)^{|G|}{\eta_1}(G'''){\eta_1}(F')
+G'''F''\big)+\\[4pt]&&\frac{(\lambda+1)(2\lambda+1)}{6}((-1)^{|G|}{\eta_1}(G^{(4)}){\eta_1}(F)+2G^{(4)}F')-\\[4pt]
&&\frac{\lambda(\lambda+1)(2\lambda+1)}{15}G^{(5)}F
\Big)\alpha_1^{\lambda+4}\hfill\text{ if }\lambda=\frac{-7\pm\sqrt{33}}{4}.
 \end{array}
\end{equation*}
\end{theorem}

To prove Theorem \ref{th1} we need the following Lemmas

\begin{lemma}\cite{bbbbk}\label{sa}
The $1$-cocycle $\Upsilon$ of $\mathcal{K}(1)$ is a coboundary if
and only if its restriction $\Upsilon'$ to $\mathfrak{vect}(1)$ is
a coboundary.
\end{lemma}

\begin{lemma}
\label{decom} As a $\mathcal{K}(n-1)-$module, we have
\begin{equation*}(\mathbb{D}^{n}_{\lambda,\mu})_{\bar 0}\simeq \mathbb{D}^{n-1}_{\lambda,\mu}
\oplus \mathbb{D}^{n-1}_{\lambda+\frac{1}{2},\mu+\frac{1}{2}}\\ \hbox{ and }\\
(\mathbb{D}^n_{\lambda,\mu})_{\bar1}\simeq\Pi(\mathbb{D}^{n-1}_{\lambda+\frac{1}{2},\mu}
\oplus \mathbb{D}^{n-1}_{\lambda,\mu+\frac{1}{2}}).
\end{equation*}
\end{lemma}
\begin{proofname}. Observe that the $\frak{ aff}(n-1|1)-$isomorphism (\ref{varphi})
 is also a $\mathcal{K}(n-1)-$isomorphism. Thus, by isomorphism (\ref{varphi}), we deduce a
$\mathcal{K}(n-1)-$isomorphism,
\begin{equation}\label{Phi}
\begin{array}{lll}\Phi_{\lambda,\mu}:&\mathbb{D}_{\lambda,\mu}^{n-1}\oplus
\mathbb{D}_{\lambda+\half,\mu+\half}^{n-1}\oplus
\Pi\left(\mathbb{D}_{\lambda,\mu+\half}^{n-1}\oplus
\mathbb{D}_{\lambda+\half,\mu}^{n-1}\right) \rightarrow
\mathbb{D}_{\lambda,\mu}^n\\[2pt]&A\longmapsto\varphi_{\mu}^{-1}\circ
A\circ\varphi_{\lambda}.
\end{array}
\end{equation}Here, we identify
the $\mathcal{K}(n-1)-$modules via the following isomorphisms:
\begin{gather*}\begin{array}{llllllll}
\Pi\left(\mathbb{D}^{n-1}_{\lambda,\mu+\frac{1}{2}}\right)&\rightarrow&
\mathrm{Hom_{diff}}\left(\mathbb{F}^{n-1}_\lambda,\Pi(\mathbb{F}^{n-1}_{\mu+\frac{1}{2}})\right)
\quad &\Pi(A)&\mapsto&\Pi\circ A,\\[10pt]
\Pi\left(\mathbb{D}^{n-1}_{\lambda+\frac{1}{2},\mu}\right)&\rightarrow&
\mathrm{Hom_{diff}}\left(\Pi(\mathbb{F}^{n-1}_{\lambda+\frac{1}{2}}),\mathbb{F}^{n-1}_{\mu}\right)
\quad &\Pi(A)&\mapsto& A\circ\Pi,\\[10pt]
\mathbb{D}^{n-1}_{\lambda+\frac{1}{2},\mu+\frac{1}{2}}&\rightarrow&
\mathrm{Hom_{diff}}\left(\Pi(\mathbb{F}^{n-1}_{\lambda+\frac{1}{2}}),\Pi(\mathbb{F}^{n-1}_{\mu+\frac{1}{2}})\right)
\quad &A&\mapsto&\Pi\circ A\circ\Pi.\\[10pt]
\end{array}
\end{gather*}
\end{proofname}

\begin{proofname} of Theorem \ref{th1}.
According to Lemma \ref{decom}, we see that $\mathrm{H}_{\rm
diff}^1(\mathfrak{vect}(1),\frak{
aff}(1);\mathbb{D}^1_{\lambda,\mu})$ can be
deduced from the spaces $\mathrm{H}_{\rm
diff}^1(\mathfrak{vect}(1),\frak{
aff}(1); \mathrm{D}_{\lambda,\mu})$:
\begin{equation}\label{cohomo}\begin{array}{ll}\mathrm{H}_{\rm
diff}^1\left(\mathfrak{vect}(1),\frak{
aff}(1);\mathbb{D}^1_{\lambda,\mu}\right)&\simeq\mathrm{H}_{\rm
diff}^1\left(\mathfrak{vect}(1),\frak{
aff}(1);
\mathrm{D}_{\lambda,\mu}\right)\oplus\mathrm{H}_{\rm
diff}^1\left(\mathfrak{vect}(1),\frak{
aff}(1);
\mathrm{D}_{\lambda+\frac{1}{2},\mu+\frac{1}{2}}\right)\oplus\\[8pt]
&\mathrm{H}_{\rm diff}^1\left(\mathfrak{vect}(1),\frak{
aff}(1);\Pi(
\mathrm{D}_{\lambda,\mu+\frac{1}{2}})\right) \oplus\mathrm{H}_{\rm
diff}^1\left(\mathfrak{vect}(1),\frak{
aff}(1);\Pi(
\mathrm{D}_{\lambda+\frac{1}{2},\mu})\right).
\end{array}\end{equation}
Hence, if $2(\mu-\lambda)\notin\{3,\,\dots,\,13\}$,
then the cohomology space ${\mathrm
H}^1_{\mathrm{diff}}(\mathcal{K}(1),\frak{
aff}(1|1);\mathbb{D}^1_{\lambda,\mu})$
vanishes. Indeed, let $\Upsilon$ be any element of
$Z_{\mathrm{diff}}^1(\mathcal{K}(1),\frak{
aff}(1|1);\mathbb{D}^1_{\lambda,\mu})$. Then by (\ref{CohoSpace2}) and
(\ref{cohomo}), up to a coboundary, the restriction of $\Upsilon$ to
$\mathfrak{vect}(1)$ vanishes, so $\Upsilon=0$ by Lemma \ref{sa}.
\medskip

For $2(\mu-\lambda)\in\{3,\,\dots,\,13\}$, we study the
$\frak{ aff}(1|1)-$invariant bilinear operators $\frak{J}_{\mu-\lambda+1}^{-1,\lambda}$.
To study any operators $\frak{J}_{\mu-\lambda+1}^{-1,\lambda}$ satisfying
$\delta(\frak{J}_{\mu-\lambda+1}^{-1,\lambda})=0$, we consider the
two components of its restriction to $\mathfrak{vect}(1)$ which we
compare with  $C_{\lambda,\mu}$ and
$C_{\lambda+{1\over2},\mu+{1\over2}}$ or $C_{\lambda+{1\over2},\mu}$
and $C_{\lambda,\mu+{1\over2}}$ depending on whether $\lambda-\mu$
is integer or semi-integer. For example, we show that
$\frak{J}_{5\over2}^{-1,\lambda}$ is a $1-$cocycle. Moreover, it is
non-trivial for $\lambda\neq-{1\over2}$ since, for
$g,\,f\in\mathbb{K}[x]$, we have
$\frak{J}_{5\over2}^{-1,\lambda}(X_g)(f)=-\theta
C_{\lambda,\lambda+2}(g,f)$.
\end{proofname}

\subsubsection{The space $\mathrm{H}^1_{\rm
diff}(\mathcal{K}(2),\mathfrak{aff}(2|1);\mathbb{D}^2_{\lambda,\mu})$}

The main result of this subsection is the following

\begin{theorem}
\label{th22}
\begin{equation}
\mathrm{H^1_{diff}}({\mathcal
K}(2),\mathfrak{aff}(2|1);\mathbb{D}^2_{\lambda,\mu})\simeq\left\{
\begin{array}{llllll}
\mathbb{R}&\text{if}\quad \mu-\lambda=1,\\[2pt]
\mathbb{R}^2&\text{if}\quad \mu-\lambda=2,\\[2pt]
0&\text{otherwise}.
\end{array}
\right.
\end{equation}
It is spanned by the following $1$-cocycles:

\begin{equation*}
  \begin{array}{llll}
\Upsilon^2_{\lambda,\lambda+1}(X_G)(F\alpha_2^{\lambda})&=&\left\{\begin{array}{ll}\eta_1\eta_2(G')F\alpha_2^{\lambda+1}
&\hfill\text{
if }\lambda\neq-{1\over2},\\[4pt]\Big(\eta_1\eta_2(G')F+
 (-1)^{|G|}\sum_{i=1}^2(-1)^{i}{\eta}_{3-i}(G')\eta_i(F)\Big)\alpha_2^{\lambda+1}&\hfill\text{
if }\lambda=-{1\over2}.
 \end{array}\right.\\[8pt]
 \Upsilon^2_{\lambda,\lambda+2}(X_G)(F\alpha_2^{\lambda})&=&\Big((2\lambda+1)\left(\frac{2\lambda}{3}G'''F-(-1)^{|G|}
 \sum_{i=1}^2{\eta}_i(G'')\eta_{i}(F)\right)
-2 \eta_2\eta_1(G')\eta_2\eta_1(F)\Big)\alpha_2^{\lambda+2}.
 \end{array}
\end{equation*}

\begin{equation*}
  \begin{array}{lll}
\widetilde{\Upsilon}^2_{\lambda,\lambda+2}(X_G)(F\alpha_2^{\lambda})&=&\left\{\begin{array}{ll}
\Big((-1)^{|G|}\sum_{i=1}^2(-1)^{i}{\eta}_{3-i}(G'')\eta_i(F)
+2\lambda\eta_2\eta_1(G'')F-\\[8pt]2\eta_2\eta_1(G')F'\Big)\alpha_2^{\lambda+2}&\hfill\text{if
}\lambda\neq-1,\\[10pt]
\Big((-1)^{|G|}\sum_{i=1}^2(-1)^{i}{\eta}_{3-i}(G')\eta_i(F')-\eta_2\eta_1(G')F'+\\[6pt]
(-1)^{|G|}\sum_{i=1}^2(-1)^{i}{\eta}_{3-i}(G'')\eta_i(F)-G''\eta_2\eta_1(F)\Big)\alpha_2^{\lambda+2}&\hfill\text{if
}\lambda=-1.\end{array}\right.
 \end{array}
\end{equation*}
\end{theorem}

To prove the Theorem above, we need first the following Proposition

\begin{prop} \cite{N}
\label{pro}
For $\lambda=0$ or $\lambda\neq\mu,$ any element of
$Z_{\rm diff}^1(\mathcal{K}(2),\mathbb{D}^2_{\lambda,\mu})$ is a
coboundary over $\mathcal{K}(2)$ if and only if at least one of its
restrictions to the subalgebras $\mathcal{K}(1)$ or $\mathcal{K}(1)^1$ is a
coboundary.
\end{prop}

\medskip

\begin{proofname} of Theorem \ref{th22}.
Note that, by Lemma \ref{aff}, the $\frak{ aff}(2|1)-$relative
$1$-cocycles are $\frak{ aff}(2|1)-$invariant bilinear differential operators
and by Proposition \ref{pro}, they
are related to the $\mathfrak{aff}(1|1)-$ relative $1$-cocycles.

\medskip

Now, let us study the relationship between any $1$-cocycle of
${\mathcal K}(2)$ and its restriction to the subalgebra ${\mathcal
K}(1)$. By Lemma \ref{decom}, we see that $\mathrm{H}_{\rm
diff}^1({\mathcal K}(1),\mathfrak{aff}(1|1);
\mathbb{D}^{2}_{\lambda,\mu})$ can be deduced from the spaces
$\mathrm{H}_{\rm diff}^1({\mathcal K}(1),\mathfrak{aff}(1|1);
\mathbb{D}^1_{\lambda,\mu})$:
{\small\begin{equation}\label{coho}\begin{array}{ll}\mathrm{H}_{\rm
diff}^1\left(\mathfrak{\mathcal
K}(1),\mathfrak{aff}(1|1);\mathbb{D}^{2}_{\lambda,\mu}\right)\simeq&\mathrm{H}_{\rm
diff}^1\left({\mathcal K}(1),\mathfrak{aff}(1|1);
\mathbb{D}^1_{\lambda,\mu}\right)\oplus\mathrm{H}_{\rm
diff}^1\left(\mathfrak{\mathcal K}(1),\mathfrak{aff}(1|1);
\mathbb{D}^1_{\lambda+\frac{1}{2},\mu+\frac{1}{2}}\right)\oplus\\[8pt]
&\mathrm{H}_{\rm diff}^1\left(\mathfrak{\mathcal
K}(1),\mathfrak{aff}(1|1);\Pi(
\mathbb{D}^1_{\lambda,\mu+\frac{1}{2}})\right)
\oplus\mathrm{H}_{\rm diff}^1\left(\mathfrak{\mathcal
K}(1),\mathfrak{aff}(1|1);\Pi(
\mathbb{D}^1_{\lambda+\frac{1}{2},\mu})\right).
\end{array}\end{equation}}
Hence, for $\lambda=0$ or $\lambda\neq\mu,$ if
$2(\mu-\lambda)\notin\{2,\,\dots,\,9\}$,  the corresponding
cohomology ${\mathrm H}^1_{\mathrm{diff}}(\mathcal{K}(2),\frak{
aff}(2|1);\mathbb{D}^2_{\lambda,\mu})$ vanish. Indeed, let
$\Upsilon$ be any element of
$Z_{\mathrm{diff}}^1(\mathcal{K}(2),\frak{
aff}(2|1);\mathbb{D}^2_{\lambda,\mu})$. Then by (\ref{CohoSpace3}) and
(\ref{coho}), up to a coboundary, the restriction of $\Upsilon$ to
$\mathfrak{aff}(1|1)$ vanishes, so $\Upsilon=0$ by Proposition
\ref{pro}.

For $2(\mu-\lambda)\in\{2,\,\dots,\,9\}$ or $\mu=\lambda\neq0$, we
study the operators $\mathbb{J}_{\mu-\lambda+1}^{-1,\lambda,2}$. To study these operators $\mathbb{J}_{\mu-\lambda+1}^{-1,\lambda,2}$
satisfying $\delta(\mathbb{J}_{\mu-\lambda+1}^{-1,\lambda,2})=0$, we
consider the two components of its restriction to
$\mathfrak{aff}(1|1)$ which we compare with
$\Upsilon^{1}_{\lambda,\mu}$ and
$\Upsilon^{1}_{\lambda+{1\over2},\mu+{1\over2}}$ or
$\Upsilon^{1}_{\lambda+{1\over2},\mu}$ and
$\Upsilon^{1}_{\lambda,\mu+{1\over2}}$ depending on whether
$\mu-\lambda$ is integer or semi-integer. A straightforward but long computation
leads to the result.
\end{proofname}


\end{document}